\documentclass[10pt]{article}

\usepackage{amsmath, amsthm, amssymb}
\usepackage{amsfonts}
\usepackage{latexsym}
\usepackage{graphicx}
\usepackage{comment}
\usepackage[ruled, noend, noline, linesnumbered]{algorithm2e}

\newcommand{\abs}[1]{\lvert#1\rvert}

\newtheorem{thm}{Theorem}[section]

\newtheorem{lem}[thm]{Lemma}
\theoremstyle{remark}
\newtheorem{rem}[thm]{Remark}
\theoremstyle{definition}
\usepackage[unicode=true, dvips, ps2pdf, bookmarks=true, colorlinks=false]{hyperref}

\makeatletter
\def\url@leostyle{%
  \@ifundefined{selectfont}{\def\UrlFont{\sf}}{\def\UrlFont{\small\ttfamily}}}
\makeatother
\usepackage{breakurl}
\usepackage{color}

\def\stackfrac{\genfrac{}{}{0pt}{}}

\begin{document}

\title{Root optimization of polynomials in the\\ number field sieve}
\author{Shi Bai\\
Research School of Computer Science\\
Australian National University\\
{\tt shih.bai@gmail.com}
\and
Richard P. Brent\\
Mathematical Sciences Institute\\
Australian National University\\
Australia
\and
Emmanuel Thom\'e\\
INRIA Nancy\\
Villers-l\`es-Nancy\\
France
}
\date{}
\maketitle

\abstract{The general number field sieve (GNFS) is the most efficient
algorithm known for factoring large integers. It consists of several
stages, the first one being polynomial selection. The quality of the
chosen polynomials in polynomial selection can be modelled in terms of
size and root properties. In this paper, we describe some algorithms for
selecting polynomials with very good root properties.}

\section{The general number field sieve}

The general number field sieve~\cite{tdotnfs-1993} is the most efficient
algorithm known for factoring large integers. It consists of several stages
including polynomial selection, sieving, filtering, linear algebra and
finding square roots.

Let $n$ be the integer to be factored. The number field sieve starts by
choosing two irreducible and coprime polynomials $f(x)$ and $g(x)$ over
$\mathbb{Z}$ which share a common root $m$ modulo $n$.  In practice, the
notations $F(x, y)$ and $G(x, y)$ for the homogenized polynomials
corresponding to $f$ and $g$ are often used.  We want to find many coprime
pairs $(a, b) \in \mathbb{Z}^2$ such that the polynomials values $F(a, b)$
and $G(a, b)$ are simultaneously smooth with respect to some upper bound
$B$. An integer is smooth with respect to bound $B$ (or $B$-smooth) if none
of its prime factors are larger than $B$. The lattice
sieving~\cite{tls-1993} and line sieving~\cite{fiwtnfs-1993} are commonly used to identify such
pairs $(a, b)$.  The running-time of sieving depends on the quality of the
chosen polynomials in polynomial selection, hence many polynomial pairs will
be generated and optimized in order to produce a best one.

\pagebreak[3]
This paper discusses algorithms for root optimization in polynomial
selection in the number field sieve. We mainly focus on polynomial selection
with two polynomials, one of which is a linear polynomial.

\section{Polynomial selection}

For large integers, most methods \cite{fiwtnfs-1993, kleinjung_2006,
kleinjung_2008, 1998-murphy-model, 1998-murphy-thesis} use a linear
polynomial for $g(x)$ and a quintic or sextic polynomial for $f(x)$. Let
$f(x) = \sum_{i = 0}^{d} c_i x^i$ and $g(x) = m_2 x - m_1$.  The standard
method to generate such polynomial pairs is to expand $n$ in base-$(m_1,
m_2)$ so $ n = \sum_{i=0}^d c_i m_1^i m_2^{d-i}$.

The running-time of sieving depends on the smoothness of the polynomial
values $\abs{F(a, b)}$ and $\abs{G(a, b)}$. Let $\Psi{(x,x^{1/u})}$ be 
the number of $x^{1/u}$-smooth integers below $x$ for some $u$. The 
Dickman-de Bruijn function $\rho (u)$~\cite{Granville_2008} is often 
used to estimate $\Psi{(x,x^{1/u})}$. It can be shown that
\begin{equation*}
\lim_{x \rightarrow \infty} \dfrac{\Psi{(x,x^{1/u})}}{x} = \rho (u).
\end{equation*}
The Dickman-de Bruijn function satisfies the differential equation
\begin{equation*}
u \rho'(u) + \rho(u-1) = 0, \quad \rho(u) = 1 \; \text{ for } 0 \le u \le 1.
\end{equation*}
It may be shown that $\rho$ satisfies the asymptotic estimate
\begin{equation*}
\log (\rho(u)) =  -(1+o(1)) u \log u \; \text{ as $u \rightarrow \infty$.}
\end{equation*}
For practical purposes, the frequency of smooth numbers can be
approximated by the Canfield-Erd\H os-Pomerance theorem, which can for
example be stated as follows~\cite{Hildebrand_1993}.
\begin{thm}
\label{smoothestimate}
For any fixed $\epsilon>0$, we have
\begin{equation*}
\Psi(x, x^{1/u}) = x u^{-u(1+o(1))}
\end{equation*}
as $x^{1/u}$ and $u$ tends to infinity, uniformly in the region $x \ge u^{u/(1-\epsilon)}$.
\end{thm}

It is desirable that the polynomial pair can produce many smooth integers
across the sieve region. This heuristically requires that the size of
polynomial values is small in general. In addition, one can choose an
algebraic polynomial $f(x)$ which has many roots modulo small prime powers.
Such a choice is driven by inheritance of practices which already date back
to the CFRAC era, where suitable multipliers were chosen precisely in order
to optimize this very property~\cite{MoBr75,PoWa83}.  Then the polynomial
values are likely to be divisible by small prime powers. This may increase
the smoothness chance for polynomial values. We describe some
methods~\cite{kleinjung_2006,1998-murphy-thesis} to estimate and compare the
quality of polynomials.

\subsection{Sieving test}

A sieving experiment over short intervals is a relatively accurate method to
compare polynomial pairs. It is often used to compare several polynomial
candidates in the final stage of the polynomial selection.
Ekkelkamp~\cite{Ekkelkamp_2008} also described a method for predicting the
number of relations needed in the sieving. The method conducts a short
sieving test and simulates relations based on the test results. Experiments
show that the prediction of the number of relations is within $2\%$ of the
number of relations needed in the actual factorization.

\subsection{Size property}

Let $(a, b)$ be pairs of relatively prime integers in the sieving region
$\Omega$. For the moment, we assume that a rectangular sieving region is
used where $\abs{a} \le U$ and $0 < b \le U$. We also assume that
polynomial values $\abs{F(a, b)}$ and $\abs{G(a, b)}$ behave like random
integers of similar size. 
The number of sieving reports (coprime pairs that lead to smooth polynomial
values) can be approximated by
\begin{equation*} \dfrac{6}{\pi^2} \iint\limits_{\Omega} \, \rho \left(
\dfrac{\log \abs{F(x, y)}}{\log B}\right) \rho \left( \dfrac{\log
\abs{G(x, y)}}{\log B}\right) \mathrm{d} x \, \mathrm{d} y
\end{equation*} The multiplier $6/\pi^2$ accounts for the probability of
$a, b$ being relatively prime.

Since $G$ is a linear polynomial, we may assume that $\log (\abs{G(a, b)})$ does not
vary much across the sieving region. A simplified approximation to compare polynomials
(ignoring the constant multiplier) is to compare
\begin{equation}
\iint\limits_{\Omega} \rho \left( \dfrac{\log \abs{F(x, y)}}{\log B}\right)
 \mathrm{d} x \, \mathrm{d} y.
\label{equ:estimate_rho}
\end{equation}

The base-$(m_1, m_2)$ expansion~\cite{kleinjung_2006, kleinjung_2008} gives
polynomials whose coefficients are $O(n^{1/(d+1)})$. The leading
coefficients $c_d$ and $c_{d-1}$ are much smaller than $n^{1/(d+1)}$. The
coefficient $c_{d-2}$ is slightly smaller than $n^{1/(d+1)}$. For such
polynomials, it is often better to use a skewed sieving region where the
sieving bounds for $a, b$ have ratio $s$, while keeping the area of the
sieving region $2U^2$. The sieving bounds become $\abs{a} \le U \sqrt{s}$
and $0 < b \le U/\sqrt{s}$. Each monomial in the polynomial is bounded by
$c_{i} U^d s^{i-d/2}$.

In the integral~(\ref{equ:estimate_rho}), computing $\rho$ is
time-consuming, especially if there are many candidates. We can use some
coarser approximations. Since $\rho(u)$ is a decreasing function of $u$, we
want to choose a polynomial pair such that the size of $\abs{F(a, b)}$ and
$\abs{G(a, b)}$ is small on average over all $(a, b)$. This roughly requires
that the coefficients of the polynomials are small in absolute value. We can
compare polynomials by the logarithm of a $L^2$-norm for polynomial $F(x,y)$ by
\begin{equation}
\dfrac{1}{2} \log \left(\displaystyle s^{-d}  \int_{-1}^{1} \int_{-1}^{1} F^2(xs, y)
 \, \mathrm{d} x \, \mathrm{d} y  \right).
\label{equation_modified_l2norm}
\end{equation}
where $s$ is the skewness of sieving region. Polynomials which minimize
the expression~\ref{equation_modified_l2norm} are expected to be better
than others.

\subsection{Root property}

If a polynomial $f(x)$ has many roots modulo small prime powers, the
polynomial values may behave more smoothly than random integers of about the
same size. Boender, Brent, Montgomery and Murphy~\cite{boender_1997,
1998-murphy-model, 1998-murphy-thesis, murphy_1998} described some
quantitative measures of this effect (root property).
\smallskip

Let $p$ be a fixed prime.
Let $\nu_p(x)$ denote the exponent of the largest power of $p$ dividing
the integer $x$ and
$\nu_p(0) = \infty$. Let $S$ be a set of integers. We use (the same)
notation ${\nu_p}(S)$ to denote the expected $p$-valuation of $x \in S$.
If integers in 
$S$ are random and uniformly distributed\footnote{We consider integer
random variables within a large enough bounded sample space.},
the expected $p$-valuation ${\nu_p}(S)$ is
\begin{equation*}
\nu_p(S)=\mathop{\mathrm{E}}\limits_{x\in S}[\nu_p(x)]=\mathrm{Pr}(\nu_p\ge1)+
\mathrm{Pr}(\nu_p\ge2)+\cdots=\frac1p+\frac1{p^2}+\cdots=
 \frac{1}{p-1}.  
\end{equation*}
Thus, in an informal (logarithmic) sense, an integer $s$ in $S$ contains an
expected power $p^{1/(p-1)}$.

Let now $S$ be a set of polynomial values $f(x)$. We use (the same) notation
${\nu_p}(S)$ (or ${\nu_p}(f)$) to denote the expected $p$-valuation of the
polynomial values $S$. Hensel's lemma gives conditions when a root of $f
\pmod {p^e}$ can be lifted to a root of $f \pmod {p^{e+1}}$.
\begin{lem}[Hensel's lemma]
Let $r_1$ be a root of $f(x)$ modulo an odd prime $p$.
\begin{enumerate}
\item If $r_1$ is a simple root,  $f(x) \pmod {p^e}$ has an unique root 
$r_e \equiv r_1 \pmod {p}$ for each $e > 1$.
\item If $r_e$ is a multiple root\footnote{We say that $r_e$ is a multiple root of
$f \pmod {p^e}$ if $f'(r_e) \equiv 0 \pmod p$.} of $f \pmod {p^e}$ for $e \ge 1$,
there are two possible cases. If $p^{e+1} \mid f(r_e)$, then 
$\forall \, i \in [0, p),$ $p^{e+1} \mid f(r_{e} + i \, p^e)$.
If $p^{e+1} \nmid f(r_e)$, $r_{e}$ cannot be lifted to a root modulo $p^{e+1}$.
\end{enumerate}
\end{lem}

Assume now that the integers $x$ leading to the values $f(x)\in S$ are
uniformly random.
There are two cases. First, suppose $p \nmid \Delta$, the discriminant of
$f(x)$. $p$ is an unramified prime. Then $f(x) \pmod p$ has only simple
roots. Let $n_p$ be the number of roots. The expected $p$-valuation of
polynomial values is $\nu_p(f) = n_p/(p-1)$ (apply the formula above,
using $\mathrm{Pr}(\nu_p\ge e)=n_p/p^e$).

The second case is when $p \mid \Delta$. Here one may get multiple roots.
The expected $p$-valuation may be obtained by counting
the number of lifted roots.
\medskip

In the number field sieve, we want to know the expected $p$-valuation of
homogeneous polynomial values $F(a,b)$, where $(a,b)$ is a pair of coprime
integers, and $F(x, y)$ is the homogenous polynomial corresponding to
$f(x)$. We assume in the following that $(a,b)$ is a uniformly random
pair of coprime integers.
We have \begin{equation}\nu_p(F(a, b))=\nu_p(F(\lambda a, \lambda
b))\label{eq:valuation-meaningful-projectively}\end{equation} for any integer $\lambda$ coprime to $p$. A pair of
coprime integers $(a,b)$ maps to a point $(a:b)$ on the projective
line $\mathbf{P}^1(\mathbb{F}_p)$. Because of
property~(\ref{eq:valuation-meaningful-projectively})
above, pairs for which
$\nu_p(F(a, b))>0$ correspond to the points of the zero-dimensional variety on
$\mathbf{P}^1(\mathbb{F}_p)$ defined by the polynomial $F$.

The projective line $\mathbf{P}^1(\mathbb{F}_p)$ has $p+1$ points,
consisting of $p$ affine points which can be represented as $(x:1)$ with
$x\in\mathbb{F}_p$, together with the point at infinity $(1:0)$. Among
these, the zeroes of $F$ correspond, for affine points $(x:1)$, to affine
roots $x\in\mathbb{F}_p$ of the dehomogenized polynomial $f$. The point
at infinity is a zero of $F$ if and only if the leading coefficient $c_d$
of $f$ cancels modulo $p$. If $F$ has a total of $n_p$ affine and
projective zeroes in $\mathbf{P}^1(\mathbb{F}_p)$, then $F(a,b)$ for
coprime $(a,b)$ is divisible by $p$ with probability $n_p/(p+1)$.

It is also possible to look at $(a,b)$ modulo a prime power $p^e$. Then
$(a,b)$ maps to an equivalence class $(a:b)$ on the projective line over
the ring $\mathbb{Z}/p^e\mathbb{Z}$. The $p$-valuation of $F$ at
$(a:b)\in\mathbf{P}^1(\mathbb{Z}/p^e\mathbb{Z})$ (an integer between $0$
and $e-1$, or ``$e$ or more'') conveys the information
of what happens modulo $p^e$. There are $p^e + p^{e-1}$ points in
$\mathbf{P}^1(\mathbb{Z}/p^e\mathbb{Z})$ ($p^e$ affine points of the form
$(x:1)$, while the remaining $p^{e-1}$ points at infinity are written as
$(1:py)$). A coprime pair $(a, b)$ chosen at random maps therefore to a
given point in
$\mathbf{P}^1(\mathbb{Z}/p^e\mathbb{Z})$ with probability
$1/(p^{e-1}(p+1))$.

Given an unramified $p$, let $F(x, y) \pmod p$ have $n_p$ affine and
projective roots (zeroes on $\mathbf{P}^1(\mathbb{F}_p)$). In application
of the Hensel Lemma (applied to $f$ at an affine root $x$, or to
$p^{d-1}f(\frac1{py})$ above the possible projective root), there is a
constant number $n_p$ of points
$(a:b)=\mathbf{P}^1(\mathbb{Z}/p^e\mathbb{Z})$ such that
$\nu_p(F(a:b))\ge e$, as $e$ grows. The expected $p$-valuation $\nu_p(F)$
is thus:
\begin{equation}
\label{equ_pvaluation}
\nu_p(F) = \sum_{e=1}^{\infty} \dfrac{n_p}{p^{e-1}(p+1)} = \dfrac{n_p p}{p^2-1}.
\end{equation}
For ramified $p$, simply counting the number $n_p$ of affine and
projective roots modulo $p$ is not sufficient to deduce $\nu_p(F)$. More
careful computation is needed modulo prime powers, which is addressed in
Sections~\ref{sec:faster-rs} and~\ref{subsec:stage2}.

Murphy~\cite[p. 49]{1998-murphy-thesis} defines the $\alpha (F)$ function to
compare the cumulative expected $p$-valuation of polynomial values to random
integers of similar size. $\alpha(F)$ can be considered as the logarithmic
benefit of using polynomials values compared to using random integers.
\begin{equation*}
\alpha(F) = \sum_{\stackfrac{p \le B}{p\ \text{prime}}} \left(\dfrac{1}{p-1} - \nu_p(F)\right)  \,
\log p.
\end{equation*}%
where the summand rewrites as $\left(1 - \dfrac{n_p  p}{p+1}\right) \,
\dfrac{\log p}{p -1}$ when $p$ is unramified.
In the number field sieve, $\alpha(F)$ is often
negative since we are interested in the case when $F(x, y)$ has more than
one root.

\subsection{Steps in polynomial selection}

Polynomial selection can be divided into three steps: polynomial generation,
size optimization and root optimization.

In the polynomial selection, we first generate good polynomials in terms of
the size property. Two efficient algorithms are given by Kleinjung
\cite{kleinjung_2006, kleinjung_2008}.  Once we have generated some
polynomial pairs $(f(x)=g(x)=m_2x-m_1)$ of relatively good size, the size
and root properties of these polynomials can be
further optimized using translation and rotation.

\begin{itemize}
\item Translation of $f(x)$ by $k$ gives a new polynomial $f_k(x)$ defined by
$f_k(x) = f(x + k)$. The root of $f_k(x)$ is $m_1/m_2 - k \pmod n$. The
linear polynomial $g_k(x)$ is $m_2 x - m_1 + k m_2$. Translation only
affects the size property.

\item Rotation by a polynomial $\lambda(x)$ gives a new polynomial $f_{\lambda(x)}
(x)$ defined by $f_{\lambda(x)}(x) = f(x) + \lambda(x) \, (m_2x - m_1)$. The
linear polynomial is unchanged $g_{\lambda(x)}(x) = g(x) = m_2x - m_1$. The
root is unchanged. $\lambda(x)$ is often a linear or quadratic polynomial,
depending on $n$ and the skewness of $f(x)$. Rotation can affect both size
and root properties.
\end{itemize}

Given a polynomial pair,
translation and rotation are used to find a polynomial of smaller
(skewed)
norm (cf Equation~(\ref{equation_modified_l2norm})). This is called size optimization.

Many polynomials can have comparable size after size optimization. We
produce and choose the best polynomials in terms of good $\alpha$-values.
This requires that the polynomials have many roots modulo small prime and
prime powers. This step is referred to as root optimization.

Given $f(x)$ (or $F(x, y)$), we can use polynomial rotation to find a
related polynomial $f_{\lambda(x)}(x)$  (or $F_{\lambda(x)}(x, y)$) which has
a smaller $\alpha$ but similar size. Polynomial rotation may also 
increase the size of trailing coefficients. However, if the skewness of the
polynomial is large, the size property of the polynomial may not be altered
significantly. Hence there is some room for rotation if the skewness is
large. As an indication of this, the skewed $L^\infty$ norm of $f$,
defined as $\max_i|s^{i-d/2}f_i|$, remains unchanged for $f_\lambda(x)$
as long as the trailing coefficients of $f_\lambda(x)$ do not dominate.
This is true for the polynomials generated by the algorithm~\cite{kleinjung_2008},
where the skewness for the polynomials is likely to be large.

We discuss some algorithms for root optimization in the following sections.

\section{Root sieve}
We focus on root optimization for quintic and sextic polynomials in this
chapter. Given a polynomial pair $(f, g)$, we want to find a rotated
polynomial with similar size but better root properties. We consider linear
rotations defined by $f_{u, v}(x) = f(x) + (u x + v) g(x)$. We want to
choose $(u, v)$ such that $f_{u, v}(x)$ has a small $\alpha$-value.

The straightforward way is to look at individual polynomials $f_{u, v}(x)$
for all possible $(u, v)$'s and compare their $\alpha$-values. This is
time-consuming and impractical since the permissible bounds on $U, V$ are
often huge.

Murphy~\cite[p. 84]{1998-murphy-thesis} describes a sieve-like procedure, namely
the root sieve, to find polynomials with good root properties. It is a
standard method to optimize the root property in the final stage of
polynomial selection. We describe Murphy's root sieve in
Algorithm~\ref{alg:murphy}. Let $B$ be the bound for small primes and $U, V$
be bounds for the linear rotation. The root sieve fills an array with
estimated $\alpha$-values. The $\alpha$-values are estimated from
$p$-valuation for small primes $p \le B$. Alternatively, it is sufficient to
calculate the summation of the weighted $p$-valuation $\nu_p (F) \log
p$ for the purpose of comparison. The idea of the root sieve is that, when
$r$ is a root of $f_{u, v}(x) \pmod {p^e}$, it is also a root of $f_{u+ip^e,
v+jp^e}(x) \pmod {p^e}$.

\begin{algorithm}
  \BlankLine
  \SetKwData{True}{true}
  \SetKwData{Where}{\textbf{where}}
  \SetKw{Break}{break}
  \SetKwInOut{Input}{input}
  \SetKwInOut{Output}{output}
  \Input{a polynomial pair $f, g$; integers $U$, $V$, $B$;}
  \Output{an array of approximated $\alpha$-values of dimension $U \times V$;}
  \BlankLine
  \For{$p \le B$, $p$\ prime}{
    \For{$e$ \Where $p^e \le B$}{
      \For{$x \in [0, p^e -1]$}{
        \For{$u \in [0, p^e -1]$}{
          compute $v$ in $f(x) + (ux+v)  g(x) \equiv 0 \pmod {p^e}$\;
          update $\nu_{p}(f_{u + i  p^e, v + j p^e})$ by sieving\;
        }
      }
    }
  }
  \caption{Murphy's root sieve}
  \label{alg:murphy}
\end{algorithm}

In general, the root sieve does not affect the projective roots
significantly. It is sufficient to only consider the affine roots'
contribution to the $\alpha$-value. In the end, we identify good slots
(those with small $\alpha$-values) in the sieving array. For each slot
(polynomial), we can compute a more accurate $\alpha$-value with a large
bound $B$ and re-optimize its size.

We consider the asymptotic complexity of Murphy's root sieve.
\begin{align*}
  \sum_{\stackfrac{p \le B}{p\ \text{prime}}}
  \left(\sum_{e=1}^{\lfloor\frac{\log B}{\log p} \rfloor} p^e p^e 
\left(O(1) + \dfrac{UV}{p^{2e}} \right) \right)
&= O\left(\frac{B^3}{\log B} \right)
+ UV
 \sum_{\stackfrac{p \le B}{p\ \text{prime}}} \left\lfloor\dfrac{\log B}{\log p}\right\rfloor\\
  & \approx UV \log B\int_2^B \dfrac{1}{\log^2 p } dp\\
  & =  O \left( UV \dfrac{B}{\log B} \right).
\end{align*}
We are interested in small primes and hence $B/ \log B$ is small. The
sieving bounds $U, V$ dominate the running-time $O(UV B / \log B)$.

\section{A faster root sieve}
\label{sec:faster-rs}

In the root sieve, we identify the number of roots of "rotated" polynomials
$f_{u, v} (x)$ for small primes and prime powers. In most cases, the roots
are simple, and hence their average $p$-valuation follows
Equation~(\ref{equ_pvaluation}). There is no need to count the lifted roots
for them. We describe a faster root sieve based on this idea.

We use the following facts based on Hensel's lemma. Suppose $r_1$ is a
simple root of $ f(x) \pmod p$. There exists a unique lifted root $r_{e}$ of
$f(x) \pmod{p^e}$ for each $e > 1$. In addition, each lifted root $r_{e}$ is
a simple root of $f(x) \pmod p$. For convenience, we say $r_e$ is a simple
root of $f(x) \pmod {p^e}$ if $f'(r_e) \not \equiv 0 \pmod p$.

Let $r_e$ be a simple root of a rotated polynomial $f_{u, v}(x) \pmod {p^e}$
for $e \ge 1$. It is clear that $f_{u + i p^e, v + j p^e}(x) \equiv f_{u,
v}(x) \pmod {p^e}$ for integers $i, j$. It follows that $r_e$ is also a
simple root of the rotated polynomials $f_{u + i p^e, v + j p^e}(x) \pmod
{p^e}$. Given a simple root $r_1$ of a polynomial $f_{u, v}(x) \pmod p$, the
contribution of the root $r_1$ to $\nu_p (F_{u, v})$ is $p/(p^2-1)$. We can
update the score\footnote{$\nu_p (F_{u, v}) \log p$, the contribution
of the root $r_1 \pmod p$ to $\alpha(F_{u, v})$.} for all rotated
polynomials $f_{u + i p, v + j p}(x)$ in a sieve.

If $r_{e}$ is a multiple root of $f(x) \pmod {p^e}$ for some $e \ge 1$,
there are two possible cases. If $f(r_{e}) \equiv 0 \pmod {p^{e+1}}$, then
$\forall l \in [0, p)$, $ f(r_{e} + l p^e) \equiv 0 \pmod {p^{e+1}}$. There
are $p$ lifted roots $r_{e+1}$ satisfying $r_{e+1} = (r_e+ l p^e) \equiv r_e
\pmod {p^e}$, $\forall l \in [0, p)$. In addition, the lifted roots
$r_{e+1}$ are multiple since $f'(r_{e+1}) \equiv 0 \pmod p$. On the other
hand, if $f(r_{e}) \not \equiv 0 \pmod {p^{e+1}}$, $r_{e}$ cannot be lifted
to a root modulo ${p^{e+1}}$.

Let $r_e$ be a multiple root of a rotated polynomial $f_{u, v}(x) \pmod
{p^e}$ for $e \ge 1$. It is also a multiple root for all rotated polynomials
$f_{u + i p^e, v + j p^e}(x) \pmod {p^e}$.

Let $r$ be a fixed integer modulo $p$. We discuss the case when $r$ is a
multiple root for some rotated polynomial $f_{u, v}(x) \pmod p$. We see that
$f(r) + (u r + v) g(r) \equiv 0 \pmod p$ and $f'(r) + u g(r) + (u r + v)
g'(r) \equiv 0 \pmod p$. Since $(u r +v) \equiv -f(r)/g(r) \pmod p$, we get
$u g^2(r) \equiv f(r) g'(r) - f'(r) g(r) \pmod p$.

Therefore, only $1$ in $p$ of $u$'s admits a multiple root at $r \pmod p$.
For the other $u$'s, we can compute $v$ and update the simple contribution
$p/(p^2-1)$ to slots in the sieve array. If $r$ is a multiple root of $f_{u,
v} (x) \pmod p$, we have to lift to count the lifted roots. We discuss the
details of the lifting method in the following sections. For the moment, we
describe the improved root sieve in Algorithm~\ref{alg:murphy_faster}.

\begin{algorithm}
  \BlankLine
  \SetKwData{True}{true}
  \SetKwData{Break}{break}
  \SetKwData{Where}{\textbf{where}}
  \SetKwInOut{Input}{input}
  \SetKwInOut{Output}{output}
  \Input{a polynomial pair $f, g$; integers $U$, $V$, $B$;}
  \Output{an array of approximated $\alpha$-values of dimension $U \times V$;}
  \BlankLine
  \For{$p \le B$, $p$\ prime}{
    \For{$x \in [0, p -1]$}{
      compute $\tilde{u}$ such that $\tilde{u}   g^2(x) \equiv f(x)   g'(x)
 - f'(x)   g(x) \pmod p$\;
      \For{$u \in [0, p -1] $}{
        compute $v$ such that $f(x) + u   x   g(x) + v g(x) \equiv 0 \pmod {p}$\;
        \If{$u \ne \tilde{u}$\;}{
          update $\nu_{p}(f_{u + i p, v + j p})$ in sieving\;
        }
        \Else{
          lift to count multiple roots of $f_{\bar{u}, \bar{v}}(x) \pmod {p^{e}}$ such
 that $(\bar{u}, \bar{v}) \equiv (u, v) \pmod p$, $\bar{u}, \bar{v} \le p^e$,
 $p^e \le B$ and then sieve\;
        }
      }
    }
  }
  \caption{A faster root sieve}
  \label{alg:murphy_faster}
\end{algorithm}

Let $r = x$ be fixed in Line 2. In Line 3, we compute $\tilde{u}$ such that
$r$ is a multiple root of $f_{\tilde{u}, v}(x)$ for some $v$. If $u \ne
\tilde{u}$, $r$ is a simple root for this $u$, and some $v$ which will be
computed in Line 5. If $u = \tilde{u}$, $f_{u, v}(x)$ admits $r$ as a
multiple root. We need to lift (up to degree $d$) to count the roots. The
running-time to do this is about
\begin{equation*}
  \sum_{\stackfrac{p \le B}{p\ \text{prime}}} \left( p \left( (p-1) \, \dfrac{UV}{p^{2}}  + O\left(\dfrac{UV}
{p^{2}} \right) \right) \right) =  O\left( UV \dfrac{B}{\log B}\right).
\end{equation*}
The asymptotic running-time has the same magnitude. In practice, however,
we benefit of 
not considering the prime powers. For comparison, Murphy's root 
sieve takes about $UV \sum_{p \le B} (\log B) /(\log p)$ operations, while Algorithm \ref{alg:murphy_faster} takes about $UV \sum_{p \le B} 1$ operations. Taking $B=200$ for instance. $\sum_{p \le 200} (\log 200) /(\log p) \approx 2705$ and $\sum_{p \le 200} 1 = 46.$

\section{A two-stage method}

We give a two-stage algorithm for the root optimization. The algorithm is
motivated by previous work by Gower~\cite{Gower2003}, Jason Papadopoulos
(personal communication), Stahlke and Kleinjung~\cite{Stahlke:2008}, who
suggested to consider congruence classes modulo small primes.

If the permissible rotation bounds $U, V$ are large, the root sieve can take
a long time for each polynomial. This is even more inconvenient if there are
many polynomials. We describe a faster method for root optimization based on
the following ideas.

A polynomial with only a few roots modulo small prime powers $p^e$ is less
likely to have a good $\alpha$-value. Therefore, rotated polynomials with
many (comparably) roots modulo small prime powers are first detected. A
further root sieve for larger prime powers can then be applied.

In the first stage, we find a (or some) good rotated pair $(u_0, v_0) \pmod
{p_1^{e_1} \cdots p_m^{e_m}}$ such that the polynomial $f_{{u_0}, v_0} (x)$
has many roots modulo (very) small prime powers $p_1^{e_1}$, $\cdots$,
$p_m^{e_m}$. Let $B_s$ be an upper bound for $p_m^{e_m}$. In the second stage,
we apply the root sieve in Algorithm~\ref{alg:murphy_faster} to the
polynomial $f_{{u_0}, v_0}(x)$ for larger prime powers up to some bound $B$.

\subsection{Stage 1}

Given $f(x)$, we want to find a rotated polynomial $f_{u_0, v_0}(x)$ which
has many roots modulo small primes and small prime powers. Let the prime
powers be $p_1^{e_1}, \cdots, p_m^{e_m}$. There are several ways to generate
$f_{u_0, v_0}(x)$.

First, we can root-sieve a matrix of pairs $(u, v)$ of size 
$(\prod_{i=1}^{m} p_i^{e_i})^2$ and pick up the best $(u, v)$ pair(s) 
as $(u_0, v_0)$. If the matrix is small, there is no need to restrict
the bound in the root sieve to be $B_s$. We can use the larger 
bound $B$. If the matrix is large, however, the root sieve might be 
slow. We describe a faster strategy.

We first find $m$ (or more\footnote{For each $p_i$, we can generate more
than one polynomial. In Stage $2$ we consider multi-sets of combinations.})
individual polynomials $f_{u_i, v_i, p_i}(x) \; (1 \le i \le m)$ each of
which has many roots modulo small $p_i^{e_i}$. The values $u_i$ and $v_i$
are bounded by $p_i^{e_i}$. We combine them to obtain a polynomial $f_{u_0,
v_0}(x) \pmod {\prod_{i=1}^{m} p_i^{e_i}}$ using the Chinese Remainder 
Theorem. The polynomial $f_{u_0, v_0}(x) \pmod {p_i^{k}}$ has the same 
number of roots as the individual polynomials 
$f_{u_i, v_i, p_i}(x) \pmod {p_i^{k}}$ for $1 \le k \le e_i$.
Hence the combined polynomial is likely to have many roots modulo
small prime powers $p_1^{e_1}, \cdots, p_m^{e_m}$.

\paragraph{Individual polynomials.} 
To find individual polynomials $f_{u_i, v_i, p_i} (x)$ that have many roots
modulo small prime powers $p_i^{e_i}$, we can root-sieve a square matrix $
p_i^{e_i} \times p_i^{e_i}$ and pick up the good pairs.

Alternatively, we use a lifting method together with a $p_i^2$-ary tree data
structure. This seems to be more efficient when $p_i^{2e_i}$ is large. For
each $p_i^{e_i}$, we construct a tree of height $e_i$ and record good $(u,
v)$ pairs during the lift. The lift is based on Hensel's lemma. For
convenience, we fix $f(x) \pmod p$ where $p=p_i$ and $e = e_i$ for some $i$.
We describe the method.

We create a root node. In the base case, we search for polynomials $f_{u,
v}(x) \pmod p$ ($u, v \in [0, p)$) which have many roots and record them in
the tree. There can be at most $p^2$ level-$1$ leaves for the root node. 
In practice, one can discard those leaves with fewer (comparably) roots
and only keep the best branch.

Let a level-$1$ leaf be $(u, v) \pmod p$. A simple root is uniquely lifted.
If the polynomial $f_{u, v} (x) \pmod p$ only gives rise to simple roots, we
already know the exact $p$-valuation of $f_{u, v}(x)$. In case of multiple
roots, we need to lift and record the lifted pairs. Assume that $f_{u, v}(x)
\pmod p$ has some multiple root $r_m$ and some simple root $r_s$. We want to
update the $p$-valuation for rotated polynomials
\begin{equation}
  f(x) + \left(\Big(u + \sum_{k=1}^{e-1} i_k p^k\Big) x +  \Big(v + \sum_{k=1}^{e-1}
 j_k p^k\Big)\right) g(x) \pmod {p^e}
  \label{equ:lifted_poly}
\end{equation}
where each $i_k, j_k \in [0, p)$. We give the following procedure for the lifting.
\begin{enumerate}
\item 
For a simple root $r_s$, we find out which of the rotated polynomials $f_{u
+ i p, v + j p} (x) \pmod {p^2}$ admit $r_s$ as a root.  If $f_{u + i p, v +
j p} (r_s) \equiv 0 \pmod {p^2}$ for some $i, j$, then
  \begin{equation}
    (i r_s + j)   g(r_s)  + f_{u, v} (r_s)/p \equiv 0 \pmod p.
  \label{equ:rotated_poly}
  \end{equation}
  Hence the set of $(i, j)$'s satisfies a linear congruence equation. For simple 
roots, there is no need to compute the lifted root. It is sufficient to update
 the $p$-valuation contributed by $r_s$ to polynomials $f_{u + i  p, v + j  p}(x)$.

\item 
Let $r_m$ be a multiple root of $f_{u, v}(x) \pmod p$. If a rotated
polynomial $f_{u + i p, v + j p} (x) \pmod {p^2}$ admits $r_m$ as a root for
some $(i, j)$, all the $\{r_m + l p\}$ ($0 \le l < p$) are also roots for
the polynomial. In addition, $f'_{u + i p, v + j p} (r_m + l p) \equiv 0
\pmod {p}$. We record the multiple roots $\{r_m + l p\}$ together with the
$\{(u + i p, v + j p)\}$ pairs. The procedure also works for the lift from
$p^e$ to $p^{e+1}$ for higher $e$'s.
\end{enumerate}

We consider the memory usage of the $p^2$-ary tree. If $r$ is a root of
$f_{u, v} (x) \pmod p$, Equation~(\ref{equ:rotated_poly}) shows a node $(u,
v) \pmod p$ gives $p$ lifted nodes $(u + ip, v + jp) \pmod {p^2}$ for some
$(i, j)$'s. Since $f_{u, v} (x) \pmod p$ can potentially have other roots
besides $r$, there could be at most $p^2$ pairs $(u + ip, v + jp) \pmod
{p^2}$. The procedure also needs to record the multiple roots for each node.
We are mainly interested in the bottom level leaves of the tree, those $(u,
v) \pmod {p^e}$. It is safe to delete the tree path which will not be used
anymore. Hence a depth-first lifting method can be used. In practice, the
memory usage is often smaller than a sieve array of size $p^{2e}$.

For each $p$, we find a polynomial that either has many simple roots or many
multiple roots which can be lifted further. Tiny primes $p$'s are more
likely to be ramified. Hence we are more likely to meet multiple roots for
tiny $p$.

\paragraph{CRTs.} 
For each $p$, we have generated some polynomial(s) rotated by $(u x +  v)
g(x) \pmod {p^e}$ which have comparably good expected $p$-valuation. For
convenience, we identify the rotated polynomial by pair $(u, v)$.

Stage $1$ repeats for prime powers $p_1^{e_1}, \cdots, p_m^{e_m}$. Let
$M=\prod_{i=1}^m p_i^{e_i}$. We generate the multi-sets combinations of
pairs $\{(u, v)\}$ and recover a set of $\{(u_0, v_0)\} \pmod M$. We fix
such a pair (rotated polynomial) $(u_0, v_0) \pmod M$.

The whole search space is an integral lattice of $\mathbb{Z}^2$. In Stage
$2$, we want to root-sieve on the sublattice points defined by $(u_0 +
\gamma M, v_0 + \beta M)$ where $(\gamma, \beta) \in \mathbb{Z}^2$. The
sublattice points are expected to give rotated polynomials with good root
properties, since the polynomials have many roots modulo $p_1^{e_1}, \cdots,
p_m^{e_m}$.

We often choose the $p_i$'s to be the smallest consecutive primes since they
are likely to contribute most to the $\alpha$-value. The exponents $e_i$ in
prime powers $p_i^{e_i}$ need some more inspection. If $e_i$ is too small,
the sieving range $(\gamma, \beta) \in \mathbb{Z}^2$ can be large. If $e_i$
is too large, $M$ is large and hence some polynomials which have good size
property might be omitted in the root sieve. One heuristic is to choose
$p_i^{e_i} \lesssim p_j^{e_j}$ for $i > j, \; i, j \in [1, m]$. To determine
$m$, one can choose $M$ to be comparable to the sieving bound $U$. Assume
that $M \approx U$. We can discard those $(u_0, v_0)$'s such that $u_0 > U$.
If $u_0$ is comparable to $U$, it is sufficient to use a line sieve for
constant rotations.
\begin{rem}
In the implementation, we may want to tune the parameters by trying several
sets of parameters such as various $p_i$'s and $e_i$'s. We can run a test
root sieve in short intervals. The set of parameters which generates the
 best score is then used.
\end{rem}

\subsection{Stage 2}
\label{subsec:stage2}

In Stage $2$, we apply the root sieve in Algorithm~\ref{alg:murphy_faster}
to polynomial $f_{{u_0}, v_0}(x)$, perhaps with some larger prime bound. In
the root sieve, one can reuse the code from Stage~$1$, where the updates of
$\alpha$-values can be batched. We describe the method as follows.

\paragraph{Sieve on sublattice.} Let $M = \prod_{i=1}^m p_i^{e_i}$ and $(u_0, v_0)$
 be fixed from Stage $1$. In the second stage, we do the root sieve for (larger)
 prime powers on the sublattice defined by $\{(u_0 + \gamma M, v_0 + \beta M)\}$
 where $\gamma, \beta \in \mathbb{Z}$.  Let $p$ be a prime and $r_k \pmod {p^k}$
 be a root of
\begin{equation*}
f(x) + \Big((u_0 + \gamma M)x +  (v_0 + \beta M) \Big) g(x) \pmod {p^k}
\end{equation*}
for some fixed integers $\gamma, \beta$. The sieve on the sublattice follows from
\begin{equation*}
f(r_k) + \left( \big(u_0 + M (\gamma + i  {p^k} ) \big)  r_k +  \big(v_0
 + M (\beta + j {p^k})\big) \right) g(r_k)\equiv 0  \pmod {p^k}
\end{equation*}
for integers $i, j \in \mathbb{Z}$. We consider the root sieve for a 
fixed prime $p$ in Algorithm~\ref{alg:murphy_faster}.

Let $f(x), g(x), M, u_0, v_0$ be fixed from Stage $1$. In
Algorithm~\ref{alg:murphy_faster}, we assume $u, r$ are fixed for the
moment. Let $p$ be a prime not dividing $M$. The sieve array has approximate
size $\lfloor U/M \rfloor \times \lfloor V/M \rfloor$. Each element
$(\gamma, \beta)$ in the sieve array stands for a point $(u_0 + \gamma M,
v_0 + \beta M)$ in $\mathbb{Z}^2$. We solve for $v$ in $f(r) + u r g(r) + v
g(r) \equiv 0 \pmod {p}$. Knowing $(u, v)$, we can solve for $(\gamma,
\beta)$ in $u \equiv u_0 + \gamma M \pmod p$ and $v \equiv v_0 + \beta M
\pmod p$, provided that $p \nmid M$.

For the moment, we fix integers $\gamma, \beta$. If  $r$ is a simple root,
it is sufficient to sieve $(\gamma + i p, \beta + jp)$ for various $(i,
j)$'s and update the $p$-valuation $p \log p/(p^2 - 1)$ to each slot. If $r$
is a multiple root, we can use a similar lifting procedure as in Stage $1$.
We describe the recursion to deal with multiple roots in
Algorithm~\ref{alg:recursion_multiple}.

\begin{algorithm}[htp]
  \BlankLine
  \SetKwData{True}{true}
  \SetKwData{Break}{break}
  \SetKwInOut{Input}{input}
  \SetKwInOut{Output}{output}
  \Input{a polynomial pair $f, g$; integers $U$, $V$, $B$; node $(u, v)$,
 tree height $e$, current level $k$, prime $p$;}
  \Output{updated $\alpha$-values array;}
  \BlankLine
  \For{multiple roots $r$ of $f_{u, v} (x) \pmod {p}$}{
    \For{$k < e$} {
      compute $(i, j)$'s in $(i r + j)   g(r)  + f_{u, v} (r) / p^k \equiv 0 \pmod p$\;
      create child nodes $(u + i p^{k}, v + j p^{k})$ with roots $\{r + l p^{k}\}, \; \forall l \in [0, p)$\;
      recursively call Algorithm~\ref{alg:recursion_multiple} on $(u, v)$'s leftmost child node\;
      change coordinates for current node $(u, v)$ and sieve\;
      delete current node and move to its sibling node or parent node\;
    }
  }
  \caption{Recursion for multiple roots}
  \label{alg:recursion_multiple}
\end{algorithm}

From Stage $1$, we know $u_0, v_0$. In Algorithm~\ref{alg:murphy_faster}, we
fix $u, r$ and solve for $v$. Given a multiple root $r$ of $f(x) \pmod
{p^k}$, we find pairs $(u', v')$ such that $f_{u', v'} (r) \equiv 0 \pmod
{p^{k+1}}$ where $u'\equiv u \pmod {p^k}$ and $v'\equiv v \pmod {p^k}$. We
can construct nodes representing the $(u', v')$ pairs together with their
roots. In the recursion, we compute the lifted nodes in a depth-first
manner. Once the maximum level $p^e$ is reached, we do the root sieve for
the current nodes and delete the nodes which have been sieved.

When a lifted tree node $(u', v') \pmod {p^k}$ is created, the number of
roots for $f_{u', v'}(x) \pmod {p^k}$ is known. In the root sieve, the
$\alpha$-scores can be updated in a batch for all the roots of $f_{u',
v'}(x) \pmod {p^k}$. For each node $(u', v')$, we also need to compute the
corresponding coordinates in the sieve array.

\paragraph{Primes $p$ dividing $M$.} 
We have assumed that $p$ is a prime not dividing $M$. From Stage $1$, $M$ is
a product of prime powers $p_i^{e_i}$ for $1 \le i \le m$. For accuracy, we
can also consider primes powers $p_i^{e_i'}$ with $e_i' \ne e_i$ such that
$p_i$ appears in the $M$. Let $r$ be root of $f_{u, v}(x) \pmod p$. If $r$
is a simple root, there is no need to consider any liftings. Hence we
consider polynomials $f_{u, v}(x) \pmod p$ which have a multiple root.

We fix some $p=p_i$ and $e=e_i$, which are used in Stage $1$. Let $u, v, p$
be fixed in Algorithm~\ref{alg:murphy_faster}. Let $e'$ be the exponent of
$p$ that we want to consider in Stage $2$. There are two cases depending on
$e'$.

If $e' \le e$, the points on the sublattice have equal scores contributed by
roots modulo $p^{e'}$. It is sufficient to look at the multiple roots modulo
$p^{k}$ for $k \le e'$. In Algorithm~\ref{alg:murphy_faster}, we either
sieve all slots of the array or do not sieve at all. Given $u, v, p, k$, if
$v \equiv v_0 \pmod {p^k}$ in $v \equiv v_0 + \beta M \pmod {p^k}$, we need
to sieve the whole array. This can be omitted because it will give the same
result for each polynomial and we only want to compare polynomials. If $v_0
\not \equiv v \pmod {p^k}$, no slot satisfies the equation. Therefore, it is
safe to skip the current iteration when $e' \le e$.

If $e' > e$, the rotated polynomials $(u, v) \pmod {p^k}$ for $e < k \le e'$
may have different behaviors. We describe some modifications in the lifting
procedure.  Let $u \equiv u_0 \pmod {p^k}$, $r$, $v_0$ be fixed in
Algorithm~\ref{alg:murphy_faster}. We compute $v$. We want to know which
points (polynomials) on the sieve array are equivalent to $(u, v) \pmod
{p^k}$.

For $k \le e$, the situation is similar to the case when $e' \le e$. If the
equation $v \equiv v_0 \pmod {p^k}$ is satisfied, we record the node $(u, v)
\pmod {p^k}$ for further liftings. There is no need to sieve since all slots
on the sieve array have equal scores for roots modulo $p^k$. If $v_0 \not
\equiv v \pmod {p}$ where $k=1$, we have neither to root-sieve nor record
the node.  Let $e < k \le e'$. If $(u', v') \pmod {p^k}$ satisfies $u'
\equiv u \pmod p$ and $v'\equiv v \pmod p$, we need
\begin{equation*}
v_0 + \beta M \equiv v' \pmod {p^k}.
\end{equation*}
The equation is solvable for $\beta$ only if
\begin{equation*}
v_0 \equiv v' \pmod {p^{k}}.
\end{equation*}
Hence it is safe to discard those $(u, v) \pmod p$ such that $u \equiv u_0
\pmod p$ but \\$v \not \equiv v_0 \pmod {p}$.

On the other hand, we consider some $k$ in $e  < k \le e'$. In the lifting
procedure, we record nodes without sieving until we reach the level-$(e+1)$
nodes. Starting from a node $(u, v)$ modulo $p^{e+1}$, that is $k > e$, we
want to solve the equation
\begin{equation*}
  v_0 + \beta M \equiv v \pmod {p^k}.
\end{equation*}
The depth-first lifting procedure shows that
\begin{equation*}
  v_0 \equiv v \pmod {p^e}.
\end{equation*}
Hence $\beta$ is solvable in the following equation
\begin{equation*}
  \dfrac{v_0 - v}{p^e} + \dfrac{M}{p^e} \beta \equiv 0 \pmod {p^{k-e}}
\end{equation*}
since $\gcd(M/p^e, p) = 1$. In the root sieve, we step the array by $\beta +
j p^{k-e}$ for various $j$.

\subsection{Further remarks and improvements}

Let $(U, V)$ be the rotation bounds for the polynomial. The root sieve in
Algorithm~\ref{alg:murphy_faster} runs asymptotically in time $U V B / \log
B$ (ignoring constant factors). In Stage $2$, the searching space is
restricted to a sublattice determined by $M=\prod_{i=1}^{m} p_i^{e_i}$,
where the parameters $p_i$'s depend on Stage $1$. Hence, the root sieve in
Stage $2$ runs in time about $U V B / (M^2 \log B)$.

In Stage $2$, the points not on the sublattice are discarded since compared
to points on the sublattice they have worse $p$-valuation for those $p$'s in
Stage $1$. We assumed that they were unlikely to give rise to polynomials
with good root properties. However, a polynomial could have good $\alpha(F)$
while some $p$ in $M$ gives a poor $p$-valuation. This often happens when
some $p'$-valuation of $p' \nmid M$, those ignored in Stage $1$, is
exceptionally good, and hence mitigates some poor $p$-valuation where $p
\mid M$.

Alternatively, we can use a root sieve to identify good rotations in Stage
$1$ for some small sieving bounds $(U', V')$.  Then we examine the pattern
of $p$-valuation of these polynomials and decide the congruence classes used
in Stage $2$.

We have ignored the size property of polynomials in the algorithms. We have
assumed that polynomials rotated by similar $(u, v)$'s have comparable size.
In practice, some trials are often needed to decide the sieving bounds $(U,
V)$. We give some further remarks regarding the implementation.

\paragraph{Block sieving.}  
The root sieve makes frequent memory references to the array. However, there
is only one arithmetic operation for each array element. The time spent on
retrieving memory often dominates. For instance, the root sieve may cause
cache misses if the sieve on $p$ steps over a large sieve array. A common
way to deal with cache misses is to sieve in blocks.

We partition the sieving region into multiple blocks each of whose size is
at most the cache size. In the root sieve, we attempt to keep each block in
the cache while many arithmetic operations are applied. The fragment of the
block sieving is described in Algorithm~\ref{alg:murphy_faster_order}.

\begin{algorithm}[htp]
  \BlankLine
  \SetKwData{True}{true}
  \SetKwData{Break}{break}
  \SetKwInOut{Input}{input}
  \SetKwInOut{Output}{output}
  \Input{a polynomial pair $f, g$; integers $U$, $V$, $B$;}
  \Output{an array of approximated $\alpha$-values in dimension $U \times V$;}
  \BlankLine
  \For{$x \le B$}{
    \For{each block}{
      \For{$p$ where $x < p \le B$}{
        $\cdots \cdots$
      }
    }
  }
  \caption{Block sieving}
\label{alg:murphy_faster_order}
\end{algorithm}

We have also changed the order of iterations to better facilitate the block
sieving.  This might give some benefits due to the following heuristic. In
Algorithm~\ref{alg:murphy_faster}, when $p$ is small, polynomial roots $x$
modulo $p$ are small. The number of roots $x \le p$ blocked for sieving is
also limited. Instead we block primes $p$. If $x$ is small, there are still
many $p$'s which can be blocked.

For multiple roots, we might need to sieve in steps $p^k$ for $k \ge 1$.
When $p^k$ is not too small, each block has only a few (or none) references.
In this case, we may use a sorting-based sieving procedure like the bucket
sieve~\cite{DBLP:conf/asiacrypt/AokiU04}.

\paragraph{Arithmetic.} The coefficients of the rotated polynomials are multiple
 precision numbers.  Since $p^e$ can often fit into a single precision integer,
 it is sufficient to use single precision in most parts of the algorithms.

The algorithms involve arithmetic on $p^k$ for all $k \le e$. It is
sufficient to store polynomial coefficients modulo $p^e$ and do the modulo
reduction for arithmetic modulo $p^k$. Let $D$ be a multiple precision
integer. In the algorithm, we use a single precision integer $S$ instead of
$D$ where $S = D \pmod {p^e}$.  If $x \equiv D \pmod {p^k}$ for $k \le e$,
it is clear that $x \equiv S \pmod {p^k}$. Hence we can use the $S$ in the
root optimization.

In addition, the range of possible $\alpha$-values is small. We may use short
integers to approximate the $\alpha$-values instead of storing floating
point numbers. This might save some memory.

\paragraph{Quadratic rotation.}

Sextic polynomials have been used in the factorizations of many large
integers such as RSA-768. Rotations by quadratic polynomials can be used for
sextic polynomials if the coefficients and skewness of the polynomials are
large. We have assumed that $W$ is small in $f_{w, u, v}(x)$ and we
restricted to use linear rotations in this section. If the permissible bound
for $W$ is large, we can use a similar idea to that in Stage $1$ to find
good sublattices in three variables. At the end of Stage~$1$, a set of
polynomials having good $\alpha$-values are found which are defined by
rotations of $(w_0, u_0, v_0)$'s. In Stage $2$, we root-sieve on the
sublattice $\{(w_0 + \delta M, u_0 + \gamma M, v_0 + \beta M)\}$ where
$\delta, \gamma, \beta \in \mathbb{Z}$.

\section{Conclusion}

Root optimization aims to produce polynomials that have many roots modulo
small primes and prime powers.  We gave some faster methods for root
optimization based on Hensel's lifting lemma and root sieve on congruences
classes modulo small prime powers. The algorithms described here have been
implemented and tested in practice. The implementation can be found in
CADO-NFS~\cite{cadonfs_2011}.

\bibliographystyle{abbrv}    
\bibliography{reference}
\end{document}